\documentstyle[10pt]{article}
\begin{document}
%
%
%
\newtheorem{theorem}      {Th\'eor\`eme}[section]
\newtheorem{theorem*}     {theorem}
\newtheorem{proposition}  [theorem]{Proposition}
\newtheorem{definition}   [theorem]{Definition}
\newtheorem{e-lemme}        [theorem]{Lemma}
\newtheorem{cor}   [theorem]{Corollaire}
\newtheorem{resultat}     [theorem]{R\'esultat}
\newtheorem{eexercice}    [theorem]{Exercice}
\newtheorem{rrem}    [theorem]{Remarque}
\newtheorem{pprobleme}    [theorem]{Probl\`eme}
\newtheorem{eexemple}     [theorem]{Exemple}
\newcommand{\preuve}      {\paragraph{Preuve}}
\newenvironment{probleme} {\begin{pprobleme}\rm}{\end{pprobleme}}
\newenvironment{remarque} {\begin{rremarque}\rm}{\end{rremarque}}
\newenvironment{exercice} {\begin{eexercice}\rm}{\end{eexercice}}
\newenvironment{exemple}  {\begin{eexemple}\rm}{\end{eexemple}}
%
%
\newtheorem{e-theo}      [theorem]{Theorem}
\newtheorem{theo*}     [theorem]{Theorem}
\newtheorem{e-pro}  [theorem]{Proposition}
\newtheorem{e-def}   [theorem]{Definition}
\newtheorem{e-lem}        [theorem]{Lemma}
\newtheorem{e-cor}   [theorem]{Corollary}
\newtheorem{e-resultat}     [theorem]{Result}
\newtheorem{ex}    [theorem]{Exercise}
\newtheorem{e-rem}    [theorem]{Remark}
\newtheorem{prob}    [theorem]{Problem}
\newtheorem{example}     [theorem]{Example}
\newcommand{\proof}         {\paragraph{Proof~: }}
\newcommand{\hint}          {\paragraph{Hint}}
\newcommand{\heuristicproof}{\paragraph{heuristic proof}}
\newenvironment{e-probleme} {\begin{e-pprobleme}\rm}{\end{e-pprobleme}}
\newenvironment{e-remarque} {\begin{e-rremarque}\rm}{\end{e-rremarque}}
\newenvironment{e-exercice} {\begin{e-eexercice}\rm}{\end{e-eexercice}}
\newenvironment{e-exemple}  {\begin{e-eexemple}\rm}{\end{e-eexemple}}
\newcommand{\reell}    {{{\rm I\! R}^l}}
\newcommand{\reeln}    {{{\rm I\! R}^n}}
\newcommand{\reelk}    {{{\rm I\! R}^k}}
\newcommand{\reelm}    {{{\rm I\! R}^m}}
\newcommand{\reelp}    {{{\rm I\! R}^p}}
\newcommand{\reeld}    {{{\rm I\! R}^d}}
\newcommand{\reeldd}   {{{\rm I\! R}^{d\times d}}}
\newcommand{\reelnn}   {{{\rm I\! R}^{n\times n}}}
\newcommand{\reelnd}   {{{\rm I\! R}^{n\times d}}}
\newcommand{\reeldn}   {{{\rm I\! R}^{d\times n}}}
\newcommand{\reelkd}   {{{\rm I\! R}^{k\times d}}}
\newcommand{\reelkl}   {{{\rm I\! R}^{k\times l}}}
\newcommand{\reelN}    {{{\rm I\! R}^N}}
\newcommand{\reelM}    {{{\rm I\! R}^M}}
\newcommand{\reelplus} {{{\rm I\! R}^+}}
\newcommand{\reelo}    {{{\rm I\! R}\setminus\{0\}}}
\newcommand{\reld}    {{{\rm I\! R}_d}}
\newcommand{\relplus} {{{\rm I\! R}_+}}
\newcommand{\1}        {{\bf 1}}

\newcommand{\cov}      {{\hbox{cov}}}
\newcommand{\sss}      {{\cal S}}
\newcommand{\indic}    {{{\rm I\!\! I}}}
\newcommand{\pp}       {{{\rm I\!\!\! P}}}
\newcommand{\qq}       {{{\rm I\!\!\! Q}}}
\newcommand{\ee}       {{{\rm I\! E}}}

\newcommand{\B}        {{{\rm I\! B}}}
\newcommand{\cc}       {{{\rm I\!\!\! C}}}
\newcommand{\HHH}        {{{\rm I\! H}}}
\newcommand{\N}        {{{\rm I\! N}}}
\newcommand{\R}        {{{\rm I\! R}}}
\newcommand{\D}        {{{\rm I\! D}}}
\newcommand{\Z}       {{{\rm Z\!\! Z}}}
\newcommand{\C}        {{\bf C}}	
\newcommand{\T}        {{\bf T}}	
\newcommand{\E}        {{\bf E}}	
\newcommand{\rfr}[1]    {\stackrel{\circ}{#1}}
\newcommand{\equiva}    {\displaystyle\mathop{\simeq}}
\newcommand{\eqdef}     {\stackrel{\triangle}{=}}
\newcommand{\limps}     {\mathop{\hbox{\rm lim--p.s.}}}
\newcommand{\Limsup}    {\mathop{\overline{\rm lim}}}
\newcommand{\Liminf}    {\mathop{\underline{\rm lim}}}
\newcommand{\Inf}       {\mathop{\rm Inf}}
\newcommand{\vers}      {\mathop{\;{\rightarrow}\;}}
\newcommand{\versup}    {\mathop{\;{\nearrow}\;}}
\newcommand{\versdown}  {\mathop{\;{\searrow}\;}}
\newcommand{\vvers}     {\mathop{\;{\longrightarrow}\;}}
\newcommand{\cvetroite} {\mathop{\;{\Longrightarrow}\;}}
\newcommand{\ieme}      {\hbox{i}^{\hbox{\smalltype\`eme}}}
\newcommand{\eqps}      {\, \buildrel \rm \hbox{\rm\smalltype p.s.} \over = \,}
\newcommand{\eqas}      {\,\buildrel\rm\hbox{\rm\smalltype a.s.} \over = \,}
\newcommand{\argmax}    {\hbox{{\rm Arg}}\max}
\newcommand{\argmin}    {\hbox{{\rm Arg}}\min}
\newcommand{\indep}{\perp\!\!\!\!\perp}
\newcommand{\abs}[1]{\left| #1 \right|}
\newcommand{\crochet}[2]{\langle #1 \,,\, #2 \rangle}
\newcommand{\espc}[3]   {E_{#1}\left(\left. #2 \right| #3 \right)}
\newcommand{\rang}{\hbox{rang}}
\newcommand{\rank}{\hbox{rank}}
\newcommand{\signe}{\hbox{signe}}
\newcommand{\sign}{\hbox{sign}}

\newcommand\hA{{\widehat A}}
\newcommand\hB{{\widehat B}}
\newcommand\hC{{\widehat C}}
\newcommand\hD{{\widehat D}}
\newcommand\hE{{\widehat E}}
\newcommand\hF{{\widehat F}}
\newcommand\hG{{\widehat G}}
\newcommand\hH{{\widehat H}}
\newcommand\hI{{\widehat I}}
\newcommand\hJ{{\widehat J}}
\newcommand\hK{{\widehat K}}
\newcommand\hL{{\widehat L}}
\newcommand\hM{{\widehat M}}
\newcommand\hN{{\widehat N}}
\newcommand\hO{{\widehat O}}
\newcommand\hP{{\widehat P}}
\newcommand\hQ{{\widehat Q}}
\newcommand\hR{{\widehat R}}
\newcommand\hS{{\widehat S}}
\newcommand\hTT{{\widehat T}}
\newcommand\hU{{\widehat U}}
\newcommand\hV{{\widehat V}}
\newcommand\hW{{\widehat W}}
\newcommand\hX{{\widehat X}}
\newcommand\hY{{\widehat Y}}
\newcommand\hZ{{\widehat Z}}

\newcommand\ha{{\widehat a}}
\newcommand\hb{{\widehat b}}
\newcommand\hc{{\widehat c}}
\newcommand\hd{{\widehat d}}
\newcommand\he{{\widehat e}}
\newcommand\hf{{\widehat f}}
\newcommand\hg{{\widehat g}}
\newcommand\hh{{\widehat h}}
\newcommand\hi{{\widehat i}}
\newcommand\hj{{\widehat j}}
\newcommand\hk{{\widehat k}}
\newcommand\hl{{\widehat l}}
\newcommand\hm{{\widehat m}}
\newcommand\hn{{\widehat n}}
\newcommand\ho{{\widehat o}}
\newcommand\hp{{\widehat p}}
\newcommand\hq{{\widehat q}}
\newcommand\hr{{\widehat r}}
\newcommand\hs{{\widehat s}}
\newcommand\htt{{\widehat t}}
\newcommand\hu{{\widehat u}}
\newcommand\hv{{\widehat v}}
\newcommand\hw{{\widehat w}}
\newcommand\hx{{\widehat x}}
\newcommand\hy{{\widehat y}}
\newcommand\hz{{\widehat z}}

\newcommand\tA{{\widetilde A}}
\newcommand\tB{{\widetilde B}}
\newcommand\tC{{\widetilde C}}
\newcommand\tD{{\widetilde D}}
\newcommand\tE{{\widetilde E}}
\newcommand\tF{{\widetilde F}}
\newcommand\tG{{\widetilde G}}
\newcommand\tH{{\widetilde H}}
\newcommand\tI{{\widetilde I}}
\newcommand\tJ{{\widetilde J}}
\newcommand\tK{{\widetilde K}}
\newcommand\tL{{\widetilde L}}
\newcommand\tM{{\widetilde M}}
\newcommand\tN{{\widetilde N}}
\newcommand\tOO{{\widetilde O}}
\newcommand\tP{{\widetilde P}}
\newcommand\tQ{{\widetilde Q}}
\newcommand\tR{{\widetilde R}}
\newcommand\tS{{\widetilde S}}
\newcommand\tTT{{\widetilde T}}
\newcommand\tU{{\widetilde U}}
\newcommand\tV{{\widetilde V}}
\newcommand\tW{{\widetilde W}}
\newcommand\tX{{\widetilde X}}
\newcommand\tY{{\widetilde Y}}
\newcommand\tZ{{\widetilde Z}}

\newcommand\ta{{\widetilde a}}
\newcommand\tb{{\widetilde b}}
\newcommand\tc{{\widetilde c}}
\newcommand\td{{\widetilde d}}
\newcommand\te{{\widetilde e}}
\newcommand\tf{{\widetilde f}}
\newcommand\tg{{\widetilde g}}
\newcommand\th{{\widetilde h}}
\newcommand\ti{{\widetilde i}}
\newcommand\tj{{\widetilde j}}
\newcommand\tk{{\widetilde k}}
\newcommand\tl{{\widetilde l}}
\newcommand\tm{{\widetilde m}}
\newcommand\tn{{\widetilde n}}
\newcommand\tio{{\widetilde o}}
\newcommand\tp{{\widetilde p}}
\newcommand\tq{{\widetilde q}}
\newcommand\tr{{\widetilde r}}
\newcommand\ts{{\widetilde s}}
\newcommand\tit{{\widetilde t}}
\newcommand\tu{{\widetilde u}}
\newcommand\tv{{\widetilde v}}
\newcommand\tw{{\widetilde w}}
\newcommand\tx{{\widetilde x}}
\newcommand\ty{{\widetilde y}}
\newcommand\tz{{\widetilde z}}

\newcommand\bA{{\overline A}}
\newcommand\bB{{\overline B}}
\newcommand\bC{{\overline C}}
\newcommand\bD{{\overline D}}
\newcommand\bE{{\overline E}}
\newcommand\bFF{{\overline F}}
\newcommand\bG{{\overline G}}
\newcommand\bH{{\overline H}}
\newcommand\bI{{\overline I}}
\newcommand\bJ{{\overline J}}
\newcommand\bK{{\overline K}}
\newcommand\bL{{\overline L}}
\newcommand\bM{{\overline M}}
\newcommand\bN{{\overline N}}
\newcommand\bO{{\overline O}}
\newcommand\bP{{\overline P}}
\newcommand\bQ{{\overline Q}}
\newcommand\bR{{\overline R}}
\newcommand\bS{{\overline S}}
\newcommand\bT{{\overline T}}
\newcommand\bU{{\overline U}}
\newcommand\bV{{\overline V}}
\newcommand\bW{{\overline W}}
\newcommand\bX{{\overline X}}
\newcommand\bY{{\overline Y}}
\newcommand\bZ{{\overline Z}}

\newcommand\ba{{\overline a}}
\newcommand\bb{{\overline b}}
\newcommand\bc{{\overline c}}
\newcommand\bd{{\overline d}}
\newcommand\be{{\overline e}}
\newcommand\bff{{\overline f}}
\newcommand\bg{{\overline g}}
\newcommand\bh{{\overline h}}
\newcommand\bi{{\overline i}}
\newcommand\bj{{\overline j}}
\newcommand\bk{{\overline k}}
\newcommand\bl{{\overline l}}
\newcommand\bm{{\overline m}}
\newcommand\bn{{\overline n}}
\newcommand\bo{{\overline o}}
\newcommand\bp{{\overline p}}
\newcommand\bq{{\overline q}}
\newcommand\br{{\overline r}}
\newcommand\bs{{\overline s}}
\newcommand\bt{{\overline t}}
\newcommand\bu{{\overline u}}
\newcommand\bv{{\overline v}}
\newcommand\bw{{\overline w}}
\newcommand\bx{{\overline x}}
\newcommand\by{{\overline y}}
\newcommand\bz{{\overline z}}

%
\newcommand{\AAA}{{\cal A}}
\newcommand{\BB}{{\cal B}}
\newcommand{\CC}{{\cal C}}
\newcommand{\DD}{{\cal D}}
\newcommand{\EE}{{\cal E}}
\newcommand{\FF}{{\cal F}}
\newcommand{\GG}{{\cal G}}
\newcommand{\HH}{{\cal H}}
\newcommand{\II}{{\cal I}}
\newcommand{\JJ}{{\cal J}}
\newcommand{\KK}{{\cal K}}
\newcommand{\LL}{{\cal L}}
\newcommand{\NN}{{\cal N}}
\newcommand{\MM}{{\cal M}}
\newcommand{\OO}{{\cal O}}
\newcommand{\PP}{{\cal P}}
\newcommand{\QQ}{{\cal Q}}
\newcommand{\RR}{{\cal R}}
\newcommand{\SS}{{\cal S}}
\newcommand{\TT}{{\cal T}}
\newcommand{\UU}{{\cal U}}
\newcommand{\VV}{{\cal V}}
\newcommand{\WW}{{\cal W}}
\newcommand{\XX}{{\cal X}}
\newcommand{\YY}{{\cal Y}}
\newcommand{\ZZ}{{\cal Z}}
\newcommand{\tbullet}{$\bullet$}
\newcommand{\ot}{\leftarrow}
\newcommand{\newblock}{}
\newcommand{\carre}{\hfill$\Box$}
\newcommand{\carreb}{\hfill\rule{0.25cm}{0.25cm}}
%
%
\newcommand{\dontforget}[1]
{{\mbox{}\\\noindent\rule{1cm}{2mm}\hfill don't forget : #1 \hfill\rule{1cm}{2mm}}\typeout{---------- don't forget : #1 ------------}}
\newcommand{\note}[2]
{ \noindent{\sf #1 \hfill \today}

\noindent\mbox{}\hrulefill\mbox{}
\begin{quote}\begin{quote}\sf #2\end{quote}\end{quote}
\noindent\mbox{}\hrulefill\mbox{}
\vspace{1cm}
}
\newcommand{\rond}[1]     {\stackrel{\circ}{#1}}
\newcommand{\rondf}       {\stackrel{\circ}{\FF}}
\newcommand{\point}[1]     {\stackrel{\cdot}{#1}}

\newcommand\relatif{{\rm \rlap Z\kern 3pt Z}}

\def\diagram#1{\def\normalbaselines{\baselineskip=Opt
\lineskip=10Opt\lineskiplimit=1pt}  \matrix{#1}}

{\huge Analyticit\'e des applications ${\cal CR}$  }

\bigskip 

{\bf Bernard COUPET ${}^a$, Sergey PINCHUK $^{b}$ et Alexandre SUKHOV $^{c}$

\bigskip

a et c: LATP, CNRS/ UMR n$^\circ$ 6632, CMI, Universit\'e de Provence, 39, rue
Joliot Curie, 13453 Marseille cedex 13.

b: Departement of Mathematics, Indiana University, Bloomington, Indiana, 47, USA.
NSF grant DMS 96 225 94.}

\bigskip

\bigskip
--------------------------------------------------------------------------------------
{\small 

{\bf R\'esum\'e.} {\it Nous \'etablissons des conditions suffisantes pour
l'analyticit\'e  d'une application ${\cal CR}$ lisse entre deux vari\'et\'es
analytiques r\'eelles}.}

\bigskip 

\bigskip

{\bf Analyticity of ${\cal CR}$ maps}

\bigskip 

{\small
{\bf Abstract.} {\it We establish   sufficient conditions for the analyticity 
of a smooth ${\cal CR}$ mapping between two real analytic manifolds.}}

\bigskip

-------------------------------------------------------------------------------------

\bigskip

\bigskip

{\bf Abriged English Version} 

\bigskip

In this note we establish two results giving
sufficient conditions for the analyticity of a smooth ${\cal CR}$ map between two 
real analytic manifolds in complex affine spaces of different dimensions.

Let $X \subset \cc^n$ be a real analytic hypersurface in a neighborhood of a 
 point $p\in X$,  defined by a real analytic function $\rho$ such that
  $\frac{\partial\rho}{\partial z_n}(p) \neq 0$. The Cauchy - Riemann operators on
$X$ are defined by  ${\cal L}_j =
\frac{\partial \rho}{\partial \overline z_n}\frac{\partial }{\partial \overline z_j}
- \frac{\partial \rho}{\partial\overline z_j}\frac{\partial }{\partial \overline
z_n}$, $j = 1,...n-1$. For every  $\alpha \in \N^{n-1}$, ${\cal L}^{\alpha}$
denotes the composition ${\cal L}^{\alpha} = {\cal L}_1^{\alpha_1}
...{\cal L}^{\alpha_{n-1}}_{n-1}$. Let $X' \subset \cc^{n'}$ be a real analytic set
in a neighborhood of a  point $p'\in X'$  defined by $X' = \{z' \in \cc^{n'}:
\rho'_k(z',\overline{z}') = 0, k = 1,...,d \}$ where every $\rho'_k$ is a real
analytic function. Suppose that $X$ is minimal at $p$. Consider a ${\cal C}^{\infty}$
smooth  ${\cal CR}$ map $f: X
\longrightarrow X'$ with $f(p) = p'$.  For any $z'$ in a neighborhood of 
$p'$ in $\cc^{n'}$ and $k = 1,...,d$  the function $ z \mapsto
\rho'_k(z',\overline{f}(z))$ is ${\cal C}^{\infty}$ smooth in a neighborhood of $p$
in $X$.    Consider  the function $H^{\alpha}_k : z'\mapsto ({\cal
L}^{\alpha}
\rho'_k(z',\overline{f}(z))) \vert_{z = p}$, holomorphic in a neighborhood of $p'$.
We introduce the notion  of the {\it  characteristic variety} of $f$
at $p$ which is by definition the germ of complex analytic variety in a neighborhood
of  $p'$ in $\cc^{n'}$ defined by ${\cal V}^p(f) = \{ z': H^{\alpha}_k(z') = 0, \alpha
\in \N^{n-1}, k = 1,...,d \}$. Our first result  is the following

\bigskip

{ THEOREM A.} {\it If $\dim {\cal V}^p(f) = 0$, $f$ is real
analytic in a neighborhood of $p$ on $X$.}

\bigskip

This statement generalizes several known results on the analyticity of ${\cal CR}$
mappings.

Our second result concerns the analyticity of a smooth ${\cal CR}$ map at
generic points. We use here the notion of the {\it transcendence degree}
$tr.deg_p(f)$ of a smooth ${\cal CR}$ map $f$ which measures the degree of
"non-analyticity" of the map. Let $V_p(f)$ denotes the intersection of (germs of)
complex analytic varieties through $(p,f(p))$ in $\cc^n \times \cc^{n'}$ contained
the graph of $f$. Then $tr.deg_p(f) := \dim V_p(f) - n$. We also make use of
  the notion of $(r,m)$-flateness of a real analytic
set $X'$ which means that $X'$ contains a real analytic submanifold of dimension $r$
biholomorphic to the cartesian product of the unit ball in $\cc^m$ and a real
analytic submanifold in $\cc^{n'-m}$.

\bigskip

{ THEOREM B.-} {\it Suppose that $tr.deg_z(f)$ is constant  and equal to $m$
in a neighborhood of $p$ in $X$ and that $r$ denotes the maximal rang of 
$f$ on this neighborhood. Then every neighborhood  of $p'= f(p)$ contains points
where  $X'$ is $(r,m)$-flat. In particular,  if $X'$ does not contain complex
subvarieties of positive dimension, $f$ is real analytic on an open dense subset in
$X$}.

\bigskip

--------------------------------------------------------------------------------------

\bigskip

{\bf 1. Introduction et r\'esultats.} 

\bigskip

Dans cette note, nous \'etablissons deux
r\'esultats qui donnent les conditions suffisantes pour l'analyticit\'e d'une
application ${\cal CR}$ lisse entre deux vari\'et\'es analytiques r\'eelles dans les
espaces complexes lin\'eaires de dimensions diff\'erentes.
 
Soit $X \subset \cc^n$ une hypersurface
analytique r\'eelle au voisinage d'un point $p \in X$ d\'efinie par une fonction
analytique r\'eelle $\rho$ telle que  $\frac{\partial\rho}{\partial z_n}(p) \neq 0$.
Supposons que $X$ est minimal en $p$, c'est - \`a -dire ne contient pas de germe
d'une hypersurface complexe. Les op\'erateurs de Cauchy - Riemann sur
$X$ sont d\'efinis par
${\cal L}_j =
\frac{\partial \rho}{\partial \overline z_n}\frac{\partial }{\partial \overline z_j}
- \frac{\partial \rho}{\partial\overline z_j}\frac{\partial }{\partial \overline
z_n}$, $j = 1,...n-1$. Pour tout  $\alpha \in \N^{n-1}$, ${\cal L}^{\alpha}$
d\'esigne  l'op\'erateur  ${\cal L}^{\alpha} = {\cal L}_1^{\alpha_1}
...{\cal L}^{\alpha_{n-1}}_{n-1}$. Soit  $X' \subset \cc^{n'}$ un ensemble analytique
r\'eel au voisinage d'un point $p'
\in X'$ d\'efini  par $X' = \{z' \in \cc^{n'}: \rho'_k(z',\overline{z}') = 0, k =
1,...,d \}$ o\`u chaque $\rho'_k$ est une fonction analytique r\'eelle.  Une
application
$f: X \longrightarrow \cc^{n'}$ est dite de Cauchy-Riemann (${\cal CR}$) si elle
est ${\cal C}^{\infty}$ et si ses composantes sont annul\'ees par les op\'erateurs de
Cauchy-Riemann. Nous consid\'erons une application ${\cal CR}$ $f: X
\longrightarrow X'$  telle que $f(p) = p'$. Pour  $z'$ au voisinage de
$p'$ dans $\cc^{n'}$ et $k = 1,...,d$  la fonction $ z \mapsto
\rho'_k(z',\overline{f}(z))$ \'etant de classe ${\cal C}^{\infty}$ au voisinage
de $p$ sur $X$, nous    consid\'erons  la
fonction $H^{\alpha}_k : z'\mapsto {\cal L}^{\alpha}
\rho'_k(z',\overline{f}(z)) \vert_{z = p}$, holomorphe au voisinage de $p'$. Nous
introduisons une   notion de {\it  vari\'et\'e caract\'eristique} de $f$ en
$p$ qui , par d\'efinition est le germe d'ensemble analytique complexe au voisinage
de $p'$ dans
$\cc^{n'}$ d\'efini par
${\cal V}^p(f) = \{ z': H^{\alpha}_k(z') = 0, \alpha \in \N^{n-1}, k = 1,...,d \}$.
Faisons remarquer que
$p'$ appartient \`a ${\cal V}^p(f)$ car pour tout $z \in X$ nous avons
 $\rho'_k(f(z),\overline{f}(z)) = 0$  et $f$ est ${\cal CR}$ avec $f(p) =
p'$.

Avec les notations ci-dessus, notre  premier r\'esultat  est: 

\bigskip

{ TH\'EOR\`EME A.-} {\it Si $\dim {\cal V}^{p}(f) = 0$, $f$ est
analytique r\'eelle au voisinage de $p$ sur $X$.}

\bigskip

Nous donnons quelques  situations 
o\`u la dimension de la vari\'et\'e caract\'eristique est \'egale \`a z\'ero , ce
qui assure l'analyticit\'e: (a) $X'$ est une hypersurface Levi-nond\'eg\'ener\'ee
dans
$\cc^n$ et $f: X \longrightarrow X'$ est un diff\'eomorphisme \cite{Le,Pi}; (b) $X'$
est une hypersurface de type essentiellement fini dans
$\cc^n$ et $f: X\longrightarrow X'$ est un diff\'eomorphisme \cite{BJT};
(c) $X'$ est une hypersurface de type essentiellement fini dans
$\cc^n$ et $f: X \longrightarrow X'$ est de multiplicit\'e finie \cite{BR,DF}.
(d) D'apr\`es le lemme de Hopf, les conditions de l'exemple (c)  sont
v\'erifi\'ees si $X$ et $X'$ sont des hypersurfaces pseudoconvexes  dans
$\cc^n$ sans  courbe holomorphe et $f$ est non-constante \cite{BR,DF}.

L'\'evaluation de la dimension de la vari\'et\'e
caract\'eristique dans les exemples (a)- (d) repose sur des calculs directs pour
des s\'eries formelles (voir par exemple \cite{BR}, p.495 -496); le
th\'eor\`eme A apparait comme un principe g\'en\'eral qui r\'eduit le probl\`eme
d'analyticit\'e \`a des calculs de nature alg\'ebrique.    Signalons  que les
conditions du th\'eor\`eme A sont aussi v\'erifi\'ees dans les cas  qui ne sont pas
consid\'er\'es dans des  travaux mentionn\'es ci-dessus; l'\'enonc\'e
semble aussi nouveau  dans le cas \'equidimensionnel.  Notre d\'emarche s'apparente
\`a une g\'en\'eralisation du principe  de r\'eflexion de \cite{Le,Pi}.  Nous
utilisons \'egalement des id\'ees de notre travail r\'ecent
\cite{CPS}. 

Le th\'eor\`eme A donne une condition suffisante pour l'analyticit\'e de $f$ au point
fix\'e $p$. En revanche, si la dimension de la vari\'et\'e caract\'eristique est
positive, nous pouvons seulement obtenir l'analyticit\'e en des points g\'en\'eriques.
D\'esignons par $V_p(f)$ le germe d'ensemble complexe analytique au point $(p,f(p))$
de
$\cc^n\times \cc^{n'}$ qui est l'intersection de tous les ensembles analytiques
complexes  contenant le graphe de $f$ au voisinage de
$(p,f(p))$. {\it Le degr\'e de transcendance}
$tr.deg_p(f)$ de $f$ en $p$ est par  d\'efinition la diff\'erence $\dim V_p(f) -n$.
L'application $f$ est analytique r\'eelle en $p$ si et seulement si $tr.deg_p(f) = 0$
\cite{BB}; d'autre part, $tr.deg_p(f)$ est une fonction  semicontinue
sup\'erieurement sur $X$. Cette semicontinuit\'e  du degr\'e de
transcendance implique l'existence d'un ferm\'e d'int\'erieur vide $\Sigma \subset
X$  tel que
$X \backslash \Sigma$ soit une r\'eunion finie d'ouverts sur lesquels $tr.deg_z(f)$
est constant. L'ensemble $X'$ est dit  $(r,m)$-plat en  $p'$ s'il existe
une sous - vari\'et\'e analytique r\'eelle $M \subset X'$, contenant $p'$, de
dimension $r$ et biholomorphe au produit cart\'esien $\Gamma \times D$ o\`u $D$ est un
domaine de $\cc^m$ et $\Gamma$ est une vari\'et\'e analytique r\'eelle.

Notre second r\'esultat est:

\bigskip

{ TH\'EOR\`EME B.-} {\it Supposons que $tr.deg_z(f)$ est constant et \'egal \`a $m$
sur un voisinage de $p$ dans $X$ et que $r$ d\'esigne le rang maximal de $f$ sur ce 
voisinage. Alors,  tout voisinage de $f(p)$ contient des points en lesquels  $X'$
est $(r,m)$-plat. En particulier,  si $X'$ ne contient pas de vari\'et\'e complexe
de dimension positive, $f$ est analytique r\'eelle sur un ensemble ouvert dense
de $X$}.

\bigskip

Dans le cas particulier o\`u $X$ et $X'$ sont des hypersurfaces strictement
pseudoconvexes, l'analyticit\'e de $f$ en points g\'en\'eriques a \'et\'e \'etabli
dans \cite{Fo}.

\bigskip

{\bf 2. Extension m\'eromorphe.} 

\bigskip

Notons par ${}_p{\cal C}^{\infty}(X)$
l'anneau des germes de fonctions de classe ${\cal C}^{\infty}$ au voisinage  de $p$
sur $X$ et par ${}_p{\cal C}^{\infty}_{\cal CR}(X)$ le sous-anneau de ${\cal
C}^{\infty}$ form\'e par des germes de fonctions ${\cal CR}$. Notons
\'egalement par ${}_p{\cal O}(X)$ (resp. ${}_p{\cal M}(X)$) l'anneau des
germes de fonctions holomorphes (resp. m\'eromorphes) au voisinage de $p$ sur $X$.
Pour $m \in \N$ nous disons qu'une fonction $h \in {}_p{\cal C}^{\infty}(X)$
appartient \`a
${}_p{\cal A}^{m}(X)$ si il existe $g \in ({}_p{\cal C}_{\cal CR}^{\infty}(X))^m$ et
une fonction
$H(z,\zeta,w)$, holomorphe au voisinage du point $(p,\overline{p},\overline{g(p)})$
dans $\cc^n(z)
\times\cc^{n}(\zeta) \times \cc^m(w)$ telle que $h(z) =
H(z,\overline{z},\overline{g(z)})$ pour $z \in X$. 
${}_p{\cal A}(X) : = \cup_{m \in \N} {}_p{\cal A}^m(X)$ est un sous - anneau de
${}_p{\cal C}^{\infty}(X)$ avec la graduation naturelle ${}_p{\cal A}(X)
\hookrightarrow {}_p{\cal A}^1(X) \hookrightarrow ...$ et tout op\'erateur ${\cal L}_j
:{}_p{\cal A}(X) 
\longrightarrow {}_p{\cal A}(X)$ d\' efinit une d\'erivation de ${}_p{\cal A}(X)$.

 Utilisons la notation $z = ('z,z_n)$, $'z\in\cc^{n-1}$. Notons par $\Delta_k(a,r)$
le polydisque dans $\cc^k$ de centre $a \in \cc^k$ et de rayon $r > 0$. Pour $r_1,r_2
> 0$   et $t\in \cc^{n-1}$ consid\'erons le disque lin\'eaire $d_t = \{ t \} \times
\Delta_1(p_n,r_2)$ et la famille  $\{ d_t \}_{t}$, $t \in \Delta_{n-1}('p,r_1)$.
Soit $U$ un voisinage de $p$ dans $\cc^n$; notons par $U^{+}$ (resp.
$U^{-}$)  l'ensemble $\{ z
\in U: \rho(z) > 0 \}$ (resp. $\{ z \in U: \rho(z) < 0 \}$). D'apr\`es la
minimalit\'e de $X$,  nous pouvons supposer qu'il existe une base
$\{ U_{j} \}_j$ de voisinages de
$p$ telle que pour chaque $j$, toute fonction holomorphe sur $U^{+}_j$ se prolonge
holomorphiquement sur
$U^-_j$. Par cons\'equent, toute fonction de ${}_p{\cal A}(X)$ d\'efinie sur $X \cap
U$ se prolonge comme une fonction analytique r\'eelle sur $U^{-}$ et antiholomorphe
sur tout disque $d_t \cap U^-$. Cela conduit imm\'ediatement au principe d'unicit\'e
suivant pour $\tau \in {}_p{\cal A}(X)$. Si, pour tout voisinage
$V$ de $p$ l'ensemble $\Sigma = \{ z \in X \cap V: \tau(z) = 0 \}$ est d'int\'erieur
non-vide, $\tau \equiv 0$ au voisinage de $p$ sur $X$. En particulier, ${}_p{\cal
A}(X)$ est un anneau int\`egre.

\bigskip

{ LEMME 1.-} {\it Les propri\'et\'es suivants sont v\'erifi\'ees:} 
\begin{itemize}
\item[(i)]{\it Soient $\varphi, \psi \neq 0,\tau \neq 0\in {}_p{\cal A}(X)$ et 
$\{z\in X :\psi(z) = 0 \} \subset \Sigma:= \{ z \in X: \tau(z) = 0 \}$. Supposons
que  $h =\varphi /\psi$ est de 
${\cal CR}$ sur $X \backslash \Sigma$. Alors, il existe  $\tilde{h} \in {}_p{\cal
M}(X)$ tel que $\tilde{h}\vert (X\backslash \Sigma) = h \vert (X \backslash
\Sigma)$.}
\item[(ii)] {\it Tout \'el\'ement de ${}_p{\cal C}^{\infty}_{\cal CR}(X)$,  
alg\'ebrique sur le corps  quotient $\tilde{\cal A}$ de l'anneau int\`egre
${}_p{\cal A}(X)$ appartient \`a $ {}_p{\cal O}(X)$.}
\end{itemize}
 
\bigskip

{\it D\'emonstration.-} (i) D'apr\`es le principe de r\'eflexion de Schwarz toute
fonction  continue dans $U^{-} \cup X$ et antiholomorphe dans $U^{-}$ se prolonge
dans
$U^{+}$ comme une fonction analytique r\'eelle et holomorphe sur l'intersection de
$U^+$ avec chaque disque de la famille
$\{ d_t
\}_t$. Cela implique que $\varphi$ et $\psi$ se prolongent comme fonctions analytiques
r\'eelles $\varphi^{+}$ et $\psi^{+}$ sur $U^{+}$, holomorphes sur tout disque.
Consid\'erons l'ensemble $\Sigma' \subset \cc^{n-1}$ des $t$ tels que l'intersection
$d_t \cap \Sigma$ soit de mesure lin\'eaire strictement positive. D'apr\`es le
th\'eor\`eme d'unicit\'e $(d_t \cap X) \subset \Sigma$ pour tout $t \in \Sigma'$ et
 $\Sigma'$ est un ferm\'e d'int\'erieur vide.  Fixons un point
$('a,a_n)
\in U^+
\backslash ( \Sigma' \times \cc )$ et un point $('a,b) \in 
X \backslash \Sigma$; la fonction $h$ se prolonge holomorphiquement au
voisinage $V$ de $('a,b)$ d'un c\^ot\'e de $X$. Fixons un point $('a,c) \in V$ tel que
$h$ est holomorphe sur un  polydisque $\Delta_{n}(('a,c),\delta)$ pour un $\delta > 0$
suffisamment petit. Il existe un domaine simplement connexe $D
\subset \cc$ contenant les points $c$ et $a_n$ et  tel que  la fonction $h^{+} : =
\varphi^{+}/\psi^{+}$ est
m\'eromorphe par rapport \`a $z_n$ sur  $\Delta_{n-1}('a,\delta)  \times D$. On
peut supposer qu'elle est holomorphe sur  $\Delta_{n}(('a,c),\delta)$.  Le
th\'eor\`eme d'uniformisation de Riemann et le th\'eor\`eme de Rothstein
\cite{Sh} impliquent que $h^+$ est une fonction m\'eromorphe au voisinage de $a$,
donc sur l'ensemble $U^{+} \backslash ( \Sigma'\times \cc )$ qui est un ouvert dense
dans $U^+$. Fixons maintenant un point de $U^{+} \cap ( \Sigma' \times \cc )$; on peut
supposer apr\`es un changement lin\'eaire de
coordonn\'ees qu'il coincide avec $0$ et  qu'il existe $0 < r < R$, un polydisque
$\Delta_n(q,r) \subset U^{+}\backslash (\Sigma' \times \cc)$ avec $q = ('q,q_n)$,
$'q = 0$ tels que $h^{+}$ est holomorphe sur
$\Delta_n(q,r)$ et $\Delta_{n-1}('q,r) \times \Delta_1(q_n,R)$ contient $0$. Pour tout
$w \in
\Delta_{n-1}('q,r)$ la fonction $h^{+}$ se prolonge m\'eromorphiquement sur $\{ w \}
\times \Delta_1(q_n,R)$ (car elle est holomorphe au voisinage de $q_n$ et
coincide avec un quotient de deux fonctions analytiques r\'eelles dans le disque).
 $h^{+}$ est donc m\'eromorphe sur  $U^+$ d'apr\`es le th\'eor\`eme de Rothstein.
Comme toute fonction holomorphe sur
$U^+$ se prolonge holomorphiquement sur $U$, $h^+$ se prolonge m\'eromorphiquement au
voisinage de $p$ (voir, par exemple, \cite{Sh}).

(ii) Soit $h \in {}_p{\cal C}^{\infty}_{\cal CR}(X)$ et $P = \sum_{j = 0}^{d-1} a_jx^j
+ x^d$, $a_j \in \tilde{\cal A}$ un polyn\^ome unitaire de degr\'e $d \geq 1$ minimal 
tel que
$P(h) = 0$ sur
$X$ en dehors de l'ensemble $\Sigma$ des z\'eros  des d\'enominateurs de ses
coefficients.  Appliquons les op\'erateurs de  Cauchy-Riemann ${\cal L}_j$ \`a
l'\'egalit\'e $P(h) = 0$. Comme
$P$ est minimal,  ${\cal L}_j(a_s) = 0$ sur $X \backslash \Sigma$ pout tout 
$s$ et chaque $a_s$ est m\'eromorphe au voisinage de $p$ d'apr\`es la partie (a). 
Il s'ensuit que $h$ est alg\'ebrique sur ${}_p{\cal M}(X)$ et le graphe de $h$ est
inclus dans un ensemble  analytique complexe de dimension pure $n$ au voisinage de
$(p,h(p))$ dans $\cc^{n+1}$. Alors  $h \in {}_p{\cal O}(X)$ d'apr\`es \cite{BB}.
\hfill \rule{0,25cm}{0,25cm} 

Pour la d\'emonstration du th\'eor\`eme B nous avons besoin de l'\'enonc\'e suivant:

\bigskip

{ LEMME 2.-} {\it Soit 
$h\in ({}_p{\cal C}^{\infty}_{\cal CR}(X))^k$ v\'erifiant  au moins une des
conditions suivantes :}
\begin{itemize}
\item[(i)]{\it Il existe $m \in \N$, $g \in
({}_a{\cal C}^{\infty}_{\cal CR}(X))^m$ et une fonction
$H(z,\zeta,\omega,w)$ , holomorphe au voisinage du point
$(a,\overline{a},\overline{g(a)},h(a))$ dans $\cc^n(z) \times \cc^n(\zeta) \times
\cc^m(\omega) \times \cc^k(w)$ et telle que $H(z,\overline{z},\overline{g(z)},h(z)) =
0$ au voisinage de $a$ sur $X$ et $H(z,\overline{z},\overline{g(z)},w)$ ne s'annule
pas identiqument sur $X \times \cc^k$; 
\item[(ii)] Il existe une fonction  $H(z,\zeta,\omega,w)$, holomorphe au voisinage
du point  $(a,\overline{a},h(a),\overline{h}(a))$ telle que
$H(z,\overline{z},h(z),\overline{h}(z)) = 0$ au voisinage de $a$ sur $X$ et
$H(z,\overline{z},\omega,w)$ ne s'annule pas identiquement sur $X \times
\cc^k(\omega) \times \cc^k(w)$. }
\end{itemize}
{\it Alors,  dans tout  voisinage de
$p$,   il existe un point $q \in X$ et une fonction
$G_q$ non identiquement nulle, holomorphe au voisinage de $(q,h(q))$  telle
que $G_q(z,h(z)) = 0$ au voisinage de $q$ sur $X$}. 

\bigskip 

{\it D\'emonstration.-} (i) Nous proc\'edons par r\'ecurrence sur l'entier $k$. Pour
$k = 1$,  consid\'erons la fonction $H(z,\overline{z},\overline{g(z)},w)$ comme une
s\'erie par rapport \`a $w\in\cc$ \`a coefficients dans ${}_p{\cal A}(X)$; il
existe un ouvert dense $S$ au voisinage de $p$ sur $X$ tel qu'en chaque point $q
\in S$  l'un des coefficients ne s'annule pas. On peut supposer que $q = 0$, $g(q)
= 0$, $h(q) = 0$.  D'apr\`es le th\'eor\`eme de
pr\'eparation de Weierstrass,  il existe un polyn\^ome unitaire
$P(z,\zeta,\omega)(w)$ en
$w$ avec des coefficients holomorphes au voisinage de
$0$ tel que
$P(z,\overline{z},\overline{g}(z))(h(z)) = 0$.
Le lemme 1 (ii) implique  que $h$ est holomorphe au voisinage de $0$.

Dans le cas g\'en\'eral, nous appliquons l'hypoth\`ese de r\'ecurrence et
le th\'eor\`eme de pr\'eparation de Weierstrass  et obtenons
l'\'equation
$h^d(z) +
\sum_{j=0}^{d-1}a_{j}(z,\overline{z},\overline{g(z)},h_1,...,h_{k-1})h^j_k(z) = 0$
 o\`u  les fonctions $a_j(z,\zeta,\omega,w_1,..,w_{k-1})$ sont holomorphes au
voisinage de $0$ et $d \geq 1$.   Appliquons les op\'erateurs ${\cal L}_s$
\`a cette \'equation.  Nous avons deux cas: (a) pour tout $j$ et $s$  on a 
${\cal L}_s(a_j) = 0$ au voisinage de $0$ et (b) il existe un point $t$
(arbitrairement proche vers $0$) tel qu'au voisinage de $t$ $h$ v\'erifie une
\'equation du degr\'e $d': 1 \leq d' < d$. R\'ep\'etant cet argument si n\'ecessaire,
nous pouvons  toujours se ramener au cas (a); nous consid\'erons donc seulement
ce cas.   Nous pouvons supposer qu'il existe $j_o$ tel que
$a_{j_0}(z,\overline{z},\overline{g(z)},w_1,...,w_{k-1})$ ne s'annule pas
identiquement sur $X \times \cc^{k-1}$. Soit
$a_{j_0} =\sum_I\alpha_I(z,\overline{z},\overline{g(z)}) h^I$ o\`u
$I = (i_1,...,i_{k-1},0)$. Si ${\cal L}_s(\alpha_I) = 0$ au voisinage de $0$ pour tout
$s$ et
$I$,  d'apr\`es le lemme 1 (i) chaque $\alpha_I$ est holomorphe et $h$ est
holomorphiquement d\'ependante. Si ${\cal L}_{s_0}(\alpha_{I_0})$ ne s'annule pas
identiquement pour  certain $s_0$ et $I_0$
, nous pouvons appliquer l'hypoth\`ese de
r\'ecurrence \`a $\sum_I{\cal L}_{s_0}(\alpha_I)h^I$ pour conclure.

(ii) Consid\'erons $H(z,\overline{z},h(z),\overline{h}(z))$ comme une
s\'erie enti\`ere $\sum_I \alpha_I(z,\overline{z},\overline{h(z)}) h^I(z)$ par rapport
\`a $h(z)$. Si un des coefficients ne s'annule pas   identiquement sur $X$, nous
pouvons appliquer la partie (i). Si chaque coefficient s'annule identiquement, nous
 fixons $I$ tel que $\alpha_I(z,\zeta,\omega)$ ne s'annule pas
identiquement et nous appliquons la partie (i) \`a 
l'\'equation $\overline{\alpha}_I(z,\overline{z},\overline{h(z)})= 0$. \hfill
\rule{0,25cm}{0,25cm}

\bigskip

{\bf 3. R\'esolution des \'equations analytiques.} 

\bigskip

Pour tout $m \in \N$ et toute
application $g\in ({}_p{\cal C}_{\cal CR}^{\infty}(X))^m$, d\'efinissons
${}_p{\cal A}_g^m(X)$ : une fonction $h$, d\'efinie au voisinage de
$(p,p')$ sur $X \times\cc^{n'}$, appartient \`a ${}_p{\cal A}_g^m(X)$ si il existe une
fonction
$H(z,\zeta,\omega,z')$ holomorphe au voisinage de
$(p,\overline{p},\overline{g(p)},p')$ dans $\cc^n(z) 
\times\cc^n(\zeta) \times \cc^{m}(\omega) \times \cc^{n'}(z')$ telle que $h(z,z') =
H(z,\overline{z},\overline{g(z)},z')$. Soit ${\cal H} = \{ h_j \}_{j=1}^{J} \subset
{}_p{\cal A}_g^m(X)$ une famille finie telle que $h_j(p,p') = 0$, $j = 1,...,J$.
Consid\'erons le syst\`eme d'\'equations analytiques $(*)$ : $h_j(z,f(z)) = 0, j =
1,...,J
$ o\`u les fonctions inconnues $f = (f_1,...,f_{n'})$ v\'erifient les conditions: $f
\in ({}_p{\cal C}_{\cal CR}^{\infty}(X))^{n'}$, $f(p) = p'$. Consid\'erons
\'egalement les germes d'ensembles analytiques complexes au voisinage de $p'$ dans
$\cc^{n'}$ d\'efinis par ${\cal U}({\cal H}) = \{ z' \in \cc^{n'}: h_j(p,z') = 0, j =
1,...,J \}$ et par ${\cal V}({\cal H}) = \{ z' \in
\cc^{n'}: {\cal L}^{\alpha} h_j(z,z') \vert_{z = p} = 0, \alpha \in \N^{n-1}, j =
1,...,J
\}$.

\bigskip  

{ PROPOSITION 3.-} {\it   Si $\dim {\cal V}({\cal H}) = 0$, toute
solution $f \in {}_p{\cal C}^{\infty}_{\cal CR}(X)$, $f(p) = p'$ de (*) est
analytique r\'eelle au voisinage de $p$ sur $X$.} 

\bigskip

{\it D\'emonstration.-} (a) Supposons d'abord que  $\dim {\cal U}({\cal H}) = 0$. Soit
$h_j(z,z') = H_j(z,\overline{z},\overline{g(z)},z')$; on peut supposer que $p = 0$,
$g(p) = 0$,
$p' = 0$. Consid\'erons  le syst\`eme  d'\'equations holomorphes $H_j(z,\zeta,w,z') =
0, j = 1,...,J$. D'apr\`es le th\'eor\`eme fondamental sur la description locale des
vari\'et\'es complexes analytiques,  il existe les polyn\^omes
$Q_j(z,\zeta,w)(x) = \sum_{\nu = 0}^{k_j} q_{js}(z,\zeta,w)x^s$, $q_{jk_j} \equiv 1$
, $k_j \geq 1$ \`a coefficients holomorphes au voisinage de
$0$ dans $\cc^n(z) \times \cc^n(\zeta) \times \cc^m(w)$ et tels que
les solutions de ce  syst\`eme v\'erifient  au voisinage de $0$ les \'equations  
$Q_j(z,\zeta,w)(z'_j) = 0$, $j = 1,...,n'$. Par cons\'equent, les solutions de
(*)  v\'erifient les \'equations $Q_j(z,\overline{z},\overline{g(z)})(f_j(z)) = 0$
pour $z\in X$, $j = 1,...,n'$ et  le lemme 1 (ii) s'applique.
 
(b) Soit $\dim {\cal V}({\cal H}) = 0$. Il existe
$k \in \N$ tel que ${\cal V}({\cal H})$ est d\'efini par les \'equations ${\cal
L}^{\alpha}h_j(z,z') \vert_{z = p} = 0$, $\vert \alpha \vert
\leq k$, $j = 1,...,J$. Pour $m' \in \N$ suffisamment grand, consid\'erons
l'application $g' \in ({}_p{\cal C}_{\cal CR}^{\infty}(X))^{m'}$ dont les
composantes sont $\frac{\partial^{\vert \beta \vert} g_j}{\partial z_1^{\beta_1}...
\partial z_n^{\beta_n}}$, $\vert \beta \vert \leq k$. Soient $h^{\alpha}_j : =
{\cal L}^{\alpha} h_j \in {}_p{\cal A}^{m'}_{g'}(X)$ et ${\cal H}' = \{
h_j^{\alpha} \}_{\vert \alpha \vert \leq k, j = 1,...,J}$. Comme $h^{\alpha}_j(z,f(z))
= 0$ et ${\cal V}({\cal H}) = {\cal U}({\cal H}')$, la partie (a) permet de
conclure. \hfill \rule{0,25cm}{0,25cm}

\bigskip

{\bf 4. D\'emonstrations des th\'eor\`emes.}

\bigskip 

 D\'emonstration du th\'eor\`eme  A.- Soient $h_j(z,z') =
\rho'_j(z',\overline{f(z)})$,
$g: = f$, ${\cal H} = \{ h_j \}_{j = 1,...d}$. La vari\'et\'e caract\'eristique
${\cal V}^p(f)$ coincide avec  ${\cal V}({\cal H})$ et la proposition 3  permet de
conclure.

\bigskip  

D\'emonstration du th\'eor\`eme B.-  Il existe un ensemble ouvert dense dans un 
voisinage de
$p$ sur $X$ tel que pour tout points
$a$ de cet ensemble $(a,f(a))$ est un point r\'egulier de $V_p(f)$. Soit $z' =
('z,''z)$ avec $'z \in \cc^m$. D'apr\`es le th\'eor\`eme des fonctions implicites, au
voisinage de $(a,f(a))$ on peut repr\'esenter
$V_p(f)$ dans la forme $''z = G(z,'z)$ o\`u $G$ est une fonction holomorphe. Par
cons\'equent,  l'intersection $V_p(f) \cap (X \times X')$ est donn\'ee par les
\'equations analytiques r\'eelles $R_j(z,\overline{z},'z,\overline{'z}) :=
\rho'_j(G(z,'z),\overline{G(z,'z)}) = 0$, $z \in X$, $j = 1,...,d$ et contient le
graphe de $f$. Comme
$tr.deg_z(f)$ est \'egal  \`a $m$ au voisinage de $p$, le lemme 2 (ii) implique
que
$R_j\equiv 0$ par rapport \`a $z \in X$ au voisinage de $a$ pour tout $'z$ fix\'e. 

Consid\'erons les projections canoniques $\pi: \cc^n \times \cc^{n'} \longrightarrow
\cc^n$ et $\pi': \cc^n \times \cc^{n'} \longrightarrow \cc^{n'}$ et posons $V_p(f
\vert X):= \pi^{-1}(X) \cap V_p(f)$. D'apr\`es ce qui pr\'ec\'ede $\pi'(V_p(f \vert
X))\subset X'$. Le rang de la restriction $\pi' \vert V_p(f \vert X)$ est sup\'erieur
\`a $r$ car $V_p(f \vert X)$ contient le graphe de $f$; de plus $\pi'$ est
biholomorphe sur chaque fibre de $\pi$.  D'apr\`es le th\'eor\`eme du rang il existe
une sous-vari\'et\'e $N\subset X$ analytique r\'eelle telle que la restriction $\pi':
\pi^{-1}(M) \cap V_p(f)\longrightarrow M: = \pi'(V_p(f\vert X)$ est un
diff\'eomorphisme.  L'holomorphie de $\pi'$  implique l'\'enonc\'e.

\end{document}